\documentclass[12pt]{amsart}
\usepackage{amsfonts,latexsym,rawfonts,amsmath,amssymb,amsthm,graphicx}

\numberwithin{equation}{section}
\newtheorem{theorem}{Theorem}[section]
\newtheorem{lem}[theorem]{Lemma}
\newtheorem{thm}[theorem]{Theorem}

\newtheorem{rem}[theorem]{Remark}
\newcounter{Cnumber}

\def\s{\,\,\,\,}

\def\dint{\displaystyle{\int}}
\def\mv{1.7ex}
\def\endproof{$\hfill\Box$\\}
\def\R{\mathbb{R}}

\def\C{\mathbb{C}}

\title{\bf Some remarks on Willmore surfaces embedded  in $\R^3$}
\author[Y. Li]
{Yuxiang Li
\\
{\small\it Department of Mathematical Sciences},\\
{\small\it Tsinghua University,}\\
{\small\it Beijing 100084, P.R.China.}\\
{\small\it Email: yxli@math.tsinghua.edu.cn.}}
\thanks{The author is partially supported by
NSFC Project 11131007 and NSFC Project 11371211.}

\date{}
\begin{document}
\maketitle

\begin{abstract}
Let $f:\C\rightarrow\R^3$
be complete  Willmore immersion
with $\int_{\Sigma}|A_f|^2<+\infty$. We will show that
if $f$ is the limit of an embedded surface sequence,
then $f$ is a plane. As an application, we prove that  if
$\Sigma_k$ is a sequence of closed Willmore surface
embedded in $\R^3$ with $W(\Sigma_k)<C$,
and if the conformal class of
$\Sigma_k$ converges in the moduli space, then we can find
a M\"obius transformation $\sigma_k$, such that a subsequence of
$\sigma_k(\Sigma_k)$ converges smoothly.
\end{abstract}

\section{Introduction}
Let
$f:\Sigma\rightarrow\R^3$ be an embedding. We define
the first and second fundamental form of $f$ as follows:
$$g=g_{ij}dx^i\otimes dx^j=df\otimes df,\s and\s A=A_{ij}dx^i
\otimes dx^j=-df\otimes dn.$$
Let $H=g^{ij}A_{ij}$  be the mean curvature, and
$K$ be the Causs curvature. It is well-known
that
\begin{equation}\label{equation.H}
\vec{H}=Hn=\Delta_gf.
\end{equation}
We say $f$ is minimal, if $H=0$, and Willmore if $H$ satisfies
the equation:
\begin{equation}\label{Willmore}
\Delta_gH+\frac{1}{2}(|H|^2-4K)H=0.
\end{equation}
Note that \eqref{Willmore} is the Euler-Langrange equation of
Willmore functional \cite{W}:
$$W(f)=\frac{1}{4}\int|H|^2d\mu_g.$$

Now, we let $f$ be an embedding from $\C$ into $\R^3$. We assume $f$ is complete,
noncompact, with
$\int_{\C}|A|^2<+\infty$. It is well-known that when $f$ is minimal,
$f$ must be a plane \cite{P-R}. In this  paper, we
will show that such a  result is also true when $f$ is Willmore:

\begin{thm}\label{main1}
Let $f:\C\rightarrow\R^3$ be an complete Willmore embedding.
If $\int_{\C}|A|^2<+\infty$, then $f(\C)$ is a plane.
\end{thm}

\begin{rem} In \cite{C-L},
Chen and Lamm has proved that any Willmore graph over $\R^2$ in $\R^3$
must be  plane, whenever it has finite
$\|A\|_{L^2}$.

Luo and Sun proved that if the Willmore functional of the Willmore graph is finite, then
$\|A\|_{L^2}$ is finite \cite{L-S}. However, this is not true for an embedded Willmore surface. For example,
helicoids are embedded minimal surfaces ($W=0$), but have infinite  $\|A\|_{L^2}$.
\end{rem}

Next, we will show that Theorem 1.1 still holds if we replace `embedding'
with `the limit of an embedding sequence':

\begin{thm}\label{main2}
Let $f:\C\rightarrow\R^3$ be a conformal complete Willmore immersion
with $\int_{\C}|A|^2<+\infty$.
If there exist $R_k\rightarrow+\infty$ and
embedding $\phi_k:D_{R_k}\rightarrow\R^3$, such that $f_k$
converges to $f$ in $C^1(D_R)$ for any $R$,
then $f(\C)$ is a plane.
\end{thm}

As an application, we will prove the following:

\begin{thm}\label{main3}
Let $\Sigma_k$ be a sequence of closed Willmore surface
embedded in $\R^3$. We assume the genus is fixed and
$W(\Sigma_k)<C$. If the conformal class of
$\Sigma_k$ is contained in a compact subset of  the moduli space, then we can find
M\"obius transformation $\sigma_k$, such that
$\sigma_k(\Sigma_k)$ converges smoothly.
\end{thm}

\begin{rem} Let $\Sigma_k$ be a Willmore surface
immersed in $\R^3$ and $\R^4$. Bernard and  Rivi\'ere
\cite{R-Y} proved that if
$$W(\Sigma_k)<\min\{8\pi,\omega_g^n\}-\delta,$$
modulo the action of the M\"obius
group, $\{\Sigma_k\}$
is compact. By results in \cite{K-L} ( see also \cite{R2}),
when $W(\Sigma_k)<\min\{8\pi,\omega_g^n\}-\delta,$
the conformal class of $\Sigma_k$ must be compact in the moduli space.
Moreover, by Li-Yau's inequality \cite{L-Y},
$\Sigma_k$ is an embedding when $W(\Sigma_k)<8\pi$.
\end{rem}

When $f$ has no branches  at $\infty$,  Theorem \ref{main1}
and Theorem \ref{main2} can be deduced directly from the
removability of singularity \cite{K-S} and the classification
of Willmore sphere in $S^3$ \cite{Br}. In fact, the results in
\cite{Br} imply the following:

\begin{lem}\label{Bryant}
Let $f:S^2\rightarrow\R^3$ be a Willmore immersion.  If
$f$ has no transversal self-intersectiones,
then $f$ is an embedding and $f(S^2)$ is a round sphere.
\end{lem}

Then, to get Theorem \ref{main1} and \ref{main2}, we only
need to prove $f$ has no branches at $\infty$. For this
sake, we will prove the following:
\begin{lem}\label{embedding.plain}
Let $f:\C\setminus D_R\rightarrow \R^3$ be a smooth conformal complete embedding with
$$\|A\|_{L^2(\C\setminus D_R)}<+\infty,\s
\varlimsup_{|z|\rightarrow +\infty}|f(z)|\cdot|A(z)|<+\infty.$$
Then
$$\theta^2\big(f(\mu_g \llcorner \C\setminus D_R),\infty\big)=1.$$
\end{lem}

\begin{lem}\label{embedding.plain2}
Let $f:\C\setminus D_R\rightarrow \R^3$ be a smooth conformal
complete immersion with
$$\|A\|_{L^2(\C\setminus D_R)}<+\infty,\s
\varlimsup_{|z|\rightarrow +\infty}|f(z)|\cdot|A(z)|<+\infty.$$
If there exists embedding
$\phi_k:D_{R_k}\setminus D_R\rightarrow \R^3$, which
converges to $f_0$ in $C^2(D_{R'}\setminus
D_R)$  for any $R'>R$,
then
$$\theta^2\big(f(\mu_g \llcorner \C\setminus D_R),\infty\big)=1.$$
\end{lem}

{\bf Acknowledgements.} The author would like to thank
Prof. Xiang Ma and Prof. Peng Wang for stimulating discussions.

\section{Complete Willmore embedding of $\C$ in $\R^3$
with $\int|A|^2<+\infty$}
In this section, we will prove Lemma \ref{embedding.plain},
and use it to prove Theorem \ref{main1}.

\subsection{The proof of Lemma \ref{embedding.plain}}
By Theorem 4.2.1 in \cite{M-S}, we may assume
\begin{equation}\label{MS1}
g=e^{2u}g_{euc},\s
with\s u=m\log|z|+\omega,
\end{equation}
where $m$ is a nonnegative integer and $\lim_{z\rightarrow\infty}
\omega(z)$ exists. Moreover, we have
\begin{equation}\label{|f|over|z|m}
\lim_{|z|\rightarrow+\infty}\frac{|f|}{|z|^{m+1}}=
\frac{e^{w(\infty)}}{m+1}.
\end{equation}
Also by Theorem 4.2.1 in \cite{M-S}, we can obtain
\begin{equation}\label{MS2}
\theta^2\big(f(\mu_g \llcorner \C\setminus D_R),\infty\big)=m+1.
\end{equation}

Let $f_k(z)=\frac{f(r_kz)}{r_k^{m+1}}$, where
$r_k\rightarrow+\infty$. Let $H_k$, $A_k$ be the mean curvature and
the second fundamental form of $f_k$ respectively. By \eqref{equation.H},
$$\Delta f_k=\frac{1}{2}\vec{H}_k|\nabla f_k|^2.$$
Since
$|\nabla f_k|=\sqrt{2}|z|^{m}e^{\omega(r_kz)}$, we have
$$\|\Delta f_k\|_{L^2(D_r\setminus D_{\frac{1}{r}})}\leq
C(r)W(f_k,D_r\setminus D_{\frac{1}{r}}),\hbox{ and }
\lim_{k\rightarrow+\infty}\left\||z|^{-m}
|\nabla f_k|-\sqrt{2}e^{\omega(\infty)}\right\|_{C^0(D_r\setminus
D_{\frac{1}{r}})}=0.$$
Noting that $W(f_k,D_r\setminus D_{\frac{1}{r}})
\leq W(f)$ and $f_k(1)\rightarrow\frac{e^{\omega(\infty)}}{m+1}$, we get
$$\|\Delta f_k\|_{L^2(D_r\setminus D_{\frac{1}{r}})}+
\|f_k\|_{W^{1,2}(D_r\setminus D_{\frac{1}{r}})}<C(r).$$
Applying
elliptic estimates, we have $\|f_k\|_{W^{2,2}
(D_r\setminus D_{\frac{1}{r}})}<C(r)$.
Thus we may assume $f_k$ converges to $f_0$ weakly in
$W^{2,2}(D_r\setminus D_{\frac{1}{r}})$.
Then  we may assume
$df_k\otimes df_k$ converges to $df_0\otimes df_0$ in $L^q(D_r\setminus
D_{\frac{1}{r}})$ for any $q>0$. Noting that
$$df_k\otimes df_k=|z|^{2m}e^{2\omega(r_kz)}g_{euc},$$
we get
$$df_0\otimes df_0=|z|^{2m}e^{2\omega(\infty)}g_{euc}.$$
Let $A_0$ be the second fundamental form of $f_0$.
Obviously,
$$\int_{D_r\setminus D_\frac{1}{r}}|A_0|^2
\leq \lim_{k\rightarrow+\infty}\int_{\C\setminus
D_{\frac{r_k}{r}}}|A|^2=0,$$
then  $\int_{\C}|A_0|^2=0$ and the imagine
$f_0$ is in a plane. Without loss of generality,
we may assume $w(\infty)=0$ and
  $f_0=(z^{m+1},c)$.\\

Next, we prove $m=0$ by contradiction. Assume $m>0$. By
\eqref{|f|over|z|m},
when $z\in D_r\setminus D_\frac{1}{r}$,
$$|A_{k}(z)|=r_k^{m+1}|A(r_kz)|=
\frac{r_k^{m+1}}{|f(r_kz)|}|f(r_kz)||A
(r_kz)|<C(r).$$
Then $\|\Delta f_k\|_{L^\infty(D_r\setminus D_\frac{1}{r})}<C$ and $f_k$  converges in fact
in $C^1(D_r\setminus D_\frac{1}{r})$.

If we set $f_k=(\varphi_k,f_k^3)$, then
$$\varphi_k\rightarrow z^{m+1},\s f_k^3\rightarrow c\s in\s
C^1(D_r\setminus D_\frac{1}{r}).$$
Let
$$
\Sigma_k=f_k(\C\setminus D_{\frac{R}{r_k}})\cap
\left((D_4\setminus D_\frac{1}{4})\times\R\right).
$$
and
$$F_k(x^1,x^2,x^3)=\sqrt{(x^1)^2+(x^2)^2}.$$
Then $F_k$ is $C^1$-smooth  on $\Sigma_k$ with
no critical points when $k$ is sufficiently large.

Obviously, $\{y\in\Sigma_k:F_k(y)=1\}$ consists of
compact $C^1$ smooth 1-dimensional manifolds.
Since $\varphi_k\rightarrow z^{m+1}$ and
$f_k$ is an embedding, $\{z:F_k=1\}$ has at least 2
components.
Let $\{F_k=1\}=\Gamma_1\cup\Gamma_2\cdots\cup\Gamma_{m'}$, where
$\Gamma_i$ are components of $\{F_k=1\}$ and $m'\geq 2$. Let
$\phi(\cdot,t)$ be the  flow generated by $\nabla F_k/|\nabla F_k|$
and put $\Omega_i=\phi(\Gamma_i,[-\frac{1}{2},2])$.
Then
$$\bigcup_i\Omega_i=\{2\geq F_k\geq\frac{1}{2}\},\s
and\s \Omega_i\cap \Omega_j=\emptyset.$$
That is to say that $\{2\geq F_k\geq\frac{1}{2}\}$ has at least
2 components, and on each component $\Omega_i$, we can find
$y_i$ such that
$F_k(y_i)=1$.

Let $y_i=f_k(z_i)$.
Recall that  for any fixed small $\epsilon$, when $k$ is sufficiently large,
we have
$$-\epsilon\leq |\varphi_k(z_i)|-|z_i|^{m+1}<\epsilon,\s i=1,2.$$
We may assume $z_1$, $z_2\in D_{1+\epsilon'}\setminus D_{1-\epsilon'}$ such that
$$\epsilon'\ll\frac{1}{2}, \s and\s
D_{1+\epsilon'}\setminus D_{1-\epsilon'}\subset \{z:\frac{3}{2}\geq
|\varphi_k(z)|\geq\frac{3}{4}\}.$$
Take a
curve $\gamma$ such that $\gamma([0,1])
\subset D_{1+\epsilon'}\setminus D_{1-\epsilon'}$, and
$\gamma(0)=z_1$, $\gamma(1)=z_2$. Then
$$f_k(\gamma(0))=y_1,\s f_k(\gamma(1))=y_2,\s
and\s f_k(\gamma)\subset\bigcup_i\Omega_i.$$
It is a contradiction to the fact that
$\Omega_1$ and $\Omega_2$ are different components.
\endproof

\subsection{The proof of Theorem \ref{main1}}
By a result of
Huber\cite{Hu}, we may assume $f$ to be conformal.
Without loss of generality, we assume $f(0)=0$.
We may assume
$\|A\|_{L^2(\C\setminus B_R)}<\epsilon$. Then by Theorem 2.10 in \cite{K-S},
$$r\|A\|_{L^\infty(B_{2r}\setminus B_r(0))}<C\|A\|_{L^2(B_{4r}\setminus
B_\frac{r}{2}(0))}$$
whenever $r>2R$.

Let $\Sigma$ be the image of embedding $f:\C\rightarrow\R^3$.
We deduce from Lemma \ref{embedding.plain} that
$$\lim_{R\rightarrow+\infty}\frac{\mu_\Sigma(B_R)}{\pi R^2}=1.$$
Let $y_0\notin\Sigma$ and $I(y)=\frac{y-y_0}{|y-y_0|^2}$. By Lemma 4.3 in \cite{K-S2},
$I(\Sigma)$ can be extended to  a smooth closed surface. It is easy to
check that $I(\Sigma)$ is an embedded Willmore sphere. By
Lemma \ref{Bryant}, $I(\Sigma)$
must be a round sphere, which implies that $\Sigma$ is a plane. Then
we get Theorem \ref{main1}.

\section{Compactness of a Willmore embedding sequence in
$\R^3$}
In this section, we first prove Lemma \ref{embedding.plain2}, then
prove Theorem \ref{main3}. Since the proof of Theorem \ref{main2}
is very similar to Theorem \ref{main1}, we omit it.

\subsection{The proof of Lemma \ref{embedding.plain2}}
We assume
$$g=e^{2u}g_{euc},\s
with\s u=m\log|z|+\omega,$$
where $m$ is a nonnegative integer and $\lim_{z\rightarrow\infty}
\omega(z)$ exists.
Similar to the proof of Lemma \ref{embedding.plain}, we
let $f_{0,n}(z)=\frac{f(r_nz)}{r_n^{m+1}}$, where
$r_n\rightarrow+\infty$.
we may assume
  $f_{0,n}$ converges to $(z^{m+1},c)$
in $C^1(D_\frac{1}{r}\setminus D_r)$.

Recall that $\phi_k$ converges to $f$ in $C^1$. Then,
we can find $k_n$, such that $\phi_{k_n}(r_nz)$
converges to $(z^{m+1},c)$. Then using the the arguments
similar as we prove Lemma \ref{embedding.plain}, we can finish the
proof of Lemma \ref{embedding.plain2}.

\subsection{The proof of Corollary \ref{main3}}

Let $f_k$ be conformal immersion of $(\Sigma,h_k)$ into $\R^3$,
where $h_k$ is a smooth metric with constant curvature. When
the genus of $\Sigma$  is 1, we assume $\mu(h_k)=1$. Since the conformal
structure induced by $h_k$ converges in the moduli space, we
may assume $h_k$ converges smoothly to $h_0$.
By results in \cite{K-L}, we may find M\"obius transformation
$\sigma_k$ and a finite set $\mathcal{S}$,
such that $\sigma(f_k)=1$ and $\sigma_k(f_k)$ converges
in $W^{2,2}_{loc}(\Sigma\setminus\mathcal{S},h_0)$,
where
$$\mathcal{S}=\{p\in \Sigma:\lim_{r\rightarrow 0}\varliminf_{
k\rightarrow+\infty}
\int_{B_r^{h_0}(z)}|A_{f_k}|^2\geq 8\pi\}.$$
Let $f_0$ be the limit, which is a branched $W^{2,2}$-conformal
immersion. Thus $f_0$  is continuous on $\Sigma$.

The following theorems will be useful,
see \cite{D-K}, \cite{R} for proofs respectively.

\begin{thm}\label{D.K.}
Let $g_k,g$ be smooth Riemannian metrics on a surface $M$,
such that $g_k \to g$ in $C^{s,\alpha}(M)$, where $s \in N$,
$\alpha \in (0,1)$. Then for each $p \in M$ there exist
neighborhoods $U_k, U$ and smooth conformal diffeomorphisms
$\varphi_k:D \to U_k$, such that $\vartheta_k \to \vartheta$
in $C^{s+1,\alpha}(\overline{D},M)$.
\end{thm}

\begin{thm}\label{R}
Let $f:D\rightarrow\R^n$ be a conformal immersion with
$g_f=e^{2u}g_{euc}$. Assume $f$ is Willmore. Then
there exists an $\epsilon_0>0$ and a $\lambda>0$, such that if
$$\int_D|A|^2dx<\epsilon_0,\s and\s |u|<\lambda,$$
then
$$\|\nabla^kn\|_{L^\infty(D_r)}\leq C(\epsilon_0,\lambda,r)
\|A\|_{L^2(D)},$$
where $\nabla=(\frac{\partial}{\partial x^1},\frac{\partial}
{\partial x^2})$.
\end{thm}

For simplicity, choose $\epsilon_0<8\pi-\delta$.

\begin{lem}\label{embedding.disk}
Let $f:D\setminus\{0\}\rightarrow \R^3$ be a   conformal Willmore immersion with
$$\mu_f(D)+\|A\|_{L^2(\C\setminus D_R)}<+\infty.$$
If there exist $r_k\rightarrow 0$ and embedding
$f_k:D\setminus D_{r_k}\rightarrow \R^3$, which
converges to $f$ in $C^1(D\setminus
D_r)$  for any $r<1$,
then for any sufficiently small $r$
\begin{equation*}\label{measure.0}
\theta^2\big(f(\mu_f \llcorner  D_r),0\big)=1.
\end{equation*}
\end{lem}

\proof Set $g=e^{2u}g_{euc}$. Using Proposition 4.1 in \cite{K-L}, $f\in W^{2,2}(D)$, and
$$u=m\log|z|+\omega(z),$$
where $m$ is positive integer and $\omega\in C^0(D)$.
Moreover,
we have
$$\lim_{|z|\rightarrow 0}\frac{|f(z)-f(0)|}{|z|^{m+1}}=\frac{e^{\omega(0)}}{m+1},$$
and
$$\theta^2\big(f(\mu_g \llcorner D_r),0\big)=m+1.$$
Without loss of generality, we assume $f(0)=0$.

Set $\tilde{f}=\frac{f}{|f|^2}$, and $\tilde{g}=d\tilde{f}\otimes d\tilde{f}$. Then $\tilde{g}=\frac{g}{|f|^4}$.
Let $A^0=A-\frac{1}{2}H g$, which is the traceless part of $A$.
It is well-known that
$$\int_D|A^0|^2d\mu_g=\int_D|\tilde{A}^0|^2d\mu_{\tilde{g}}.$$

Put
$\tilde{g}=e^{2\tilde{u}}g_{euc}$.
We have $\tilde{u}=u-\log|f|^2$. By Gauss curvature equation
$$-\Delta \tilde{u}=\tilde{K}e^{2\tilde{u}},$$
we get
$$\begin{array}{lll}
\dint_{D_\delta\setminus D_r}\tilde{K}e^{2\tilde{u}}
&=&-\dint_{\partial D_\delta}\frac{\partial\tilde{u}}{\partial r}
+\int_{\partial D_r}\frac{\partial\tilde{u}}{\partial r}\\[\mv]
&=&-\dint_{\partial D_\delta}\frac{\partial\tilde{u}}{\partial r}
+\int_{\partial D_r}\frac{\partial u}{\partial r}
-\int_{\partial D_r}\frac{2f_rf}{|f|^2}\\[\mv]
&=&\dint_{D_\delta\setminus D_r}Ke^{2u}-\int_{\partial D_r}\frac{2f_rf}{|f|^2}.
\end{array}$$
Since
$$
\left|\int_{\partial D_r}2\frac{f_rf}{|f|^2}\right|\leq 2
\int_{\partial D_r}\frac{|f_r|}{|f|}=2\int_0^{2\pi}
\frac{e^{u(re^{i\theta})}}{r^{m}}\frac{r^{m+1}}{|f(re^{i\theta})|}d\theta
<C,$$
we get  $|\int_{D}\tilde{K}d\mu_{\tilde{g}}|<C$.
Then  $\int_{D}|\tilde{A}|^2d\mu_{\tilde{f}}<+\infty$.
Therefore, $\tilde{f}(\frac{1}{z})$
satisfies the conditions of Lemma \ref{embedding.plain2}.

Set $\hat{f}(z)=\tilde{f}(1/z)$, and $\hat{g}=d\hat{f}\otimes
d\hat{f}=e^{2\hat{u}}g_{euc}$. We have
$$\hat{u}(z)=\tilde{u}(\frac{1}{z})-2\log|z|=
m\log|\frac{1}{z}|-\log|f(\frac{1}{z})|^2-2\log|z|.$$
Then
$$\lim_{z\rightarrow\infty}
(\hat{u}(z)-m\log|z|)=
-\lim_{z\rightarrow\infty}\log \frac{|f(\frac{1}{z})|^2}{|\frac{1}{z}|^{2m+2}}
=-2\omega(0)+2\log(m+1).$$
Applying Lemma \ref{embedding.plain2} \eqref{MS2} and \eqref{MS1}, we get
$m=0$.

\endproof

We define
$$\mathcal{S'}=\{z\in \Sigma:\lim_{r\rightarrow 0}\varliminf_{
k\rightarrow+\infty}
\int_{B_r^{h_0}(z)}|A_{f_k}|^2>\frac{\epsilon_0}{2}\}.$$
We need to prove $\mathcal{S'}$ is empty.

Assume $\mathcal{S}'$ is not empty.
Given a point $p\in\mathcal{S}'$, we choose
 $U_k, U,\vartheta_k,\vartheta$ as in Theorem \ref{D.K.}, and assume
$p=0$. We can choose $U_k$  such that
$U_k\cap\mathcal{S}'=\{p\}$. Let
$$
\hat{f}_k=f_k\circ\vartheta_k
$$
and note that
$\hat{f}_k$ is a conformal map from $D$ into $\R^3$.
Let
$$
g_{\hat{f}_k}=e^{2\hat{u}_k}g_{euc},\s h_k=e^{2v_k}g_{euc}.
$$

Note that
$0$ is the only point in $D$ which satisfies
$$\lim_{r\rightarrow 0}\varliminf_{k\rightarrow+\infty}
\int_{D_r(z)}|A_{\hat{f}_k}|d\mu_{\hat{f}_k}>\frac{\epsilon_0}{2}.$$
Put
$$\int_{D_{r_k}(z_k)}|A_{f_k}|^2=\frac{\epsilon_0}{2},\s and\s
\int_{D_{r}(z)}|A_{f_k}|^2<\frac{\epsilon_0}{2},\s \forall
D_r(z)\subset D_\frac{1}{2},\s r<r_k.$$
Then  $z_k\rightarrow 0$ and $r_k\rightarrow 0$.
Let $f_k'=\frac{\hat{f}_k(r_kz+z_k)-\hat{f}_k(z_k)}{\lambda_k}$, where
$$\lambda_k=diam(\hat{f}(z_k+[0,1/2])).$$
By Theorem \ref{convergence},
$$\|u_k'\|_{L^\infty(D_r)}
\leq C(r),\s \forall r\in(0,1).$$
Then, by  Theorem \ref{R}, $f_k'$
converges smoothly  on $D_\frac{3}{4}$.

For any point $z_0\in\partial D_\frac{1}{2}$, put
$$\gamma_k(t)=\frac{1}{4}(1+t)z_0,\s t\in[0,1],\s
and\s
\tau_k=diam(f_k'(\gamma_k)).$$
Then by Theorem \ref{convergence}
and  Theorem \ref{R},
$\frac{f_k'-f_k'(z_0)}{\tau_k}$ converges smoothly.
Since $f_k'$ converges in $D_\frac{3}{4}$,
we may assume
$\tau_k\rightarrow \tau_0>0$. Then
$f_k'$ converges smoothly on $D_\frac{3}{4}(z_0)$. Thus
$f_k'$ converges smoothly on $D_1$.
In this way, we can prove that a subsequence of $f_k'$ converges
smoothly on $D_R$ for any $R$. Let $f_0'$ be the limit. Then
$u_0'\in L^\infty_{loc}(\C)$ and
$$\int_{D}|A_{f_0'}|^2=\frac{\epsilon_0}{2}.$$
Obviously, $f_0'$ is proper.
If $diam(f_0')=+\infty$, then $f_0'$ is  noncompact and complete. Then
by Theorem \ref{main2},  $f_0'$ is a plane which implies
that $\int_D|A_{f_0'}|^2=0$. A contradiction. So, $diam(f_0')<+\infty$,
then by Simon's inequality \cite{S}, $\mu(f_0')<+\infty$.
By Proposition 4.1 in \cite{K-L}, $f_0'$ can be considered
as a continuous map from $S^2$ into $\R^3$.

Now, we set $\hat{f}_k'(z)=\hat{f}_k(z_k+z)$ and
$$\mathcal{S}(\hat{f}_k')=\{
z\in \C\setminus\{0\}:\lim_{r\rightarrow 0}
\varliminf_{k\rightarrow+\infty}\int_{D_r(z)}|A_{\hat{f}_k'}|^2>
\frac{\epsilon_0}{2}\},$$
and
$$\Gamma(\theta_1,\theta_2,t)=
\{te^{i\theta}:\theta_1\leq\theta\leq\theta_2\}.$$
Since $\mathcal{S}(\hat{f}_k')$ is a finite set, we can choose $\theta_1<\theta_2$, such that
\begin{equation}\label{noca}
\left(\cup_{t\in[\frac{r_k}{r},r]}\Gamma(\theta_1,
\theta_2,t)\right)\cap\mathcal{S}(\hat{f}_k')=\emptyset.
\end{equation}
Take $t_k\in[\frac{r_k}{r},r]$, such that
$$\lambda_k'=diam(\hat{f}_k'(\Gamma(\theta_1,\theta_2,t_k)))
=\inf_{t\in[\frac{r_k}{r},r]}diam(\hat{f}_k'(\Gamma(\theta_1,
\theta_2,t))).$$
By Proposition 4.1 in \cite{K-L},
$$\lim_{t\rightarrow 0}\lim_{k\rightarrow +\infty}diam(\hat{f}_k'(\Gamma(\theta_1,\theta_2,t)))=0,\s
\lim_{t\rightarrow \infty}\lim_{k\rightarrow+\infty}diam(f_k'(\Gamma(\theta_1,\theta_2,t)))=0.$$
Then
$$t_k\rightarrow 0,\s and \s \frac{t_k}{r_k}\rightarrow+\infty.$$
Let
$$f_k''=\frac{\hat{f}_k(t_kz+z_k)
-\hat{f}_k(t_ke^{i\theta_1}+z_k)}{\lambda_k'}$$
and
$$\mathcal{S}(\{f_k''\})=\{
z\in \C\setminus\{0\}:\lim_{r\rightarrow 0}
\varliminf_{k\rightarrow+\infty}\int_{D_r(z)}|A_{\hat{f}_k}|^2
>\frac{\epsilon_0}{2}\}.$$
By \eqref{noca} and Theorem \ref{R} and Theorem \ref{convergence}, $f_k$ converges smoothly near
$\Gamma(\theta_1,\theta_2,1)$.
Following the method we get $f_0'$,
we obtain that $f_k''$
converges smoothly on any compact  subset of
$\C\setminus(\{\mathcal{S}(f_k'')\}\cup\{0\})$.
Let $f_0''$ be the limit. Then
\begin{equation}\label{infdiam}
diam(f_0''(\Gamma(\theta_1,\theta_2,t)))
=\inf_{t\in(0,\infty)}\,diam(f_0''(\Gamma(\theta_1,\theta_2,t))).
\end{equation}
Then $\mu_{f_0''}(D_r(0))=\infty$ and
$\mu_{f_0''}(\C\setminus D_r(0))=\infty$
for any $r$. Otherwise, by Proposition 4.1 in \cite{K-L},
$$\lim_{t\rightarrow 0}diam(f_0''(\Gamma(\theta_1,\theta_2,t)))=0,\s or\s
\lim_{t\rightarrow+\infty}diam(f_0''
(\Gamma(\theta_1,\theta_2,t)))=0.$$
It contradicts  \eqref{infdiam}. Thus $f_0''$ is complete, noncompact
and has at least 2 ends.

Now, choose $y_0$ such that $$d(y_0,\frac{f_k(\Sigma)-\hat{f}_k(t_ke^{i\theta_1}+z_k)}
{\lambda_k'})>\delta>0$$
Set $I=\frac{y-y_0}{|y-y_0|^2}$. Then $I(f_k'')$
converges to  $I(f_0'')$ smoothly on any compact subset of
$\C\setminus(\{0\}\cup \mathcal{S}(\{f_k''\}))$.
For any small $r$ and any $z\in \mathcal{S}(\{f_k''\})\cup\{0\}$,
 since $I(f_k'')$ is Willmore on $D_r(z)$ and
converges smoothly to $I(f_0'')$ on $\partial D_r(z)$,  we get
$Res(I(f_0''),z)=0$
(for the definition of $Res$, one can refer to \cite{K-S2}).
Then, by  Lemma 4.1 in
\cite{K-S2} (see also Theorem I.6 in \cite{R}) and
Lemma \ref{embedding.disk}, $I(f_0'')$ is
a smooth Willmore embedding on $D_r(z)$.
Moreover, for a large $R$, since $I(f_k'')$ is Willmore on
$D_R$ and converges smoothly on $\partial D_R$,
$Res(I(f_0''),\infty)$ is also 0. Then $I(f_0'')(\frac{1}{z})$ is
a smooth Willmore embedding on $D_\frac{1}{R}$.

Therefore, $I(f_0'')$ can be considered as a smooth
conformal immersion from $S^2$ into $\R^3$.
Obviously,
$I(f_0'')$ has no transversal self-intersections.
By Lemma \ref{Bryant}, $I(f_0'')$ must be a round sphere.
It contradicts the fact that $f_0''$ has at least 2 ends.

Hence we get $\mathcal{S}'=\emptyset$.\\

Then, using the argument in \cite{K-L}, we get
$\|u_k\|_{L^\infty(\Sigma)}<C$ (this can also be deduced
from Theorem \ref{convergence}).
Given a point $p\in\Sigma$, we choose
 $U_k, U,\vartheta_k,\vartheta$ as in Theorem \ref{D.K.}, and assume
$p=0$. Let $\hat{f}_k=f_k\circ\vartheta_k$, which is
conformal.
Then we can choose an $r$, such that
$\int_{D_r}|A_{\hat{f}}|^2<\epsilon_0$.
Using Theorem \ref{R}, $\hat{f}$ converges smoothly
on $D_\frac{r}{2}$. We can choose $r$ to be sufficiently
small, such that their exists $r_p$, such that
$B^{h_0}_{r_p}(p)\subset\varphi_k(D_\frac{r}{2})$. Thus
$f_k$ converges smoothly on $B^{h_0}_{r_p}(p)$.
\endproof

\section{Appendix}
The proof of the following theorem can be found in \cite{L}. But for the
convenience of the readers, we provide a proof in this appendix.
\begin{thm}\label{convergence}
Let
$f_k:D\rightarrow\R^n$ be a smooth conformal immersion which  satisfies
\begin{itemize}
\item[{\rm 1)}] $\int_{D}|A_{f_k}|^2d\mu_{f_k}<\gamma_n-\tau$,
where $\tau>0$ and
$$\gamma_n=\left\{\begin{array}{ll}
8\pi&\mbox{ when }n=3\\[\mv]
4\pi&\mbox{ when }n\geq4.\end{array}\right.$$
\item[{\rm 2)}] $f_k(D)$ can be extended to a closed
immersed surface $\Sigma_k$ with
$$\int_{\Sigma_k}|A_{f_k}|^2d\mu_{f_k}<\Lambda.$$
\end{itemize}
Take a curve $\gamma:[0,1]\rightarrow D$,
and set $\lambda_k=diam\, f_k(\gamma)$.
Then
we can find  a subsequence of $f_k'=\frac{f_k
-f_k(\gamma(0))}{\lambda_k}$ which
converges  weakly in
$W^{2,2}_{loc}(D)$.
Let $df_k'\otimes df_k'=e^{2u_k'}(dx^1\otimes dx^1+dx^2\otimes dx^2)$.
For any $r<1$,
$$\|u_k'\|_{L^\infty(D_r)}<C(r).$$
\end{thm}

\proof Put $f_k'=\frac{f_k-f_k(\gamma(0))}{\lambda_k}$,
$\Sigma_k'=\frac{\Sigma_k-f_k(\gamma(0))}{\lambda_k}$.
We have two cases:\vspace{0.7ex}

\noindent Case 1: $diam(f_k')<C$. By
inequality (1.3) in \cite{S} with $\rho=\infty$,
$\frac{\Sigma_k'\cap B_\sigma(\gamma(0))}{\sigma^2}\leq C$ for any $\sigma>0$.
Hence we get $\mu(f_k')<C$ by taking
$\sigma=diam(f_k')$. Then by Helein's convergence theorem
\cite{H,K-L},
$f_k'$ converges weakly
in $W^{2,2}_{loc}(D)$. Since
$diam\, f_k'
(\gamma)=1$, the weak limit is not trivial.\\

\noindent Case 2: $diam(f_k')\rightarrow +\infty$. We take a point
$y_0\in\R^n$ and a constant $\delta>0$, s.t.
$$B_\delta(y_0)\cap \Sigma_k'=\emptyset,\s \forall k.$$
Let $I=\frac{y-y_0}{|y-y_0|^2}$, and
$$f_k''=I(f_k'),\s \Sigma_k''=I(\Sigma_k').$$
By conformal invariance of Willmore functional
\cite{C,W}, we have
$$\int_{\Sigma_k''}|A_{\Sigma_k''}|^2d\mu_{\Sigma_k''}
=\int_{\Sigma_k}|A_{\Sigma_k}|^2d\mu_{\Sigma_k}<\Lambda.$$
Since $\Sigma_k''\subset B_\frac{1}{\delta}(0)$, also by (1.3) in \cite{S},
we get $\mu(f_k'')<C$. Let
$$\mathcal{S}(\{f_k''\}):=
\{p\in D: \lim\limits_{r\rightarrow 0}\varliminf\limits_{k\rightarrow+\infty}
\int_{D_(p)}|A_{f_k''}|^2d\mu_{f_k''}\geq \gamma_n \}.$$
Then
$f_k''$ converges weakly in $W^{2,2}_{loc}(D\setminus
\mathcal{S}(f_k''))$.

Next, we prove that $f_k''$ does not converge to a point by assumption.
If $f_k''$ converges to a point in
$W^{2,2}_{loc}(D\setminus \mathcal{S}(f_k''))$,
then the limit must be 0,  for $diam\,(f_k')$
converges to $+\infty$.
By the
definition of $f_k''$, we can find a $\delta_0>0$,
such that $f_k''(\gamma)\cap
B_{\delta_0}(0)=\emptyset$. Thus for any $p\in \gamma([0,1])
\setminus \mathcal{S}(f_k'')$, $f_k''$ will not converge to $0$. A contradiction.

Then we only need to prove that $f_k'$ converges weakly in
$W^{2,2}_{loc}(D,\R^n)$.
Let $f_0''$ be the limit of $f_k''$ which is a branched immersion of $D$ in $\R^n$.
Let $\mathcal{S}^*=f_0^{''-1}(\{0\})$,
which is isolate. Note that for any $z_0\in\mathcal{S}^*$,
there exits $m>0$, such that
$$\lim_{|z-z_0|\rightarrow 0}
\frac{|f(z)-f(z_0)|}{|z-z_0|^m}>0.$$

First, we prove that for any $\Omega\subset\subset D\setminus
(\mathcal{S}^*\cup\mathcal{S}(\{f_k''\}))$, $f_k'$
converges weakly in $W^{2,2}(D,\R^n)$:
Since $f_0''$ is continuous
on $\bar{\Omega}$, we may assume
$dist(0,f_0''(\Omega))>\delta>0$. Then $dist(0,f_k''(\Omega))>\frac{\delta}{2}$
when $k$ is sufficiently large. Noting that $f_k'
=\frac{f_k''}{|f_k''|^2}+y_0$, we get that $f_k'$ converges weakly in
$W^{2,2}(\Omega,\R^n)$.

Next, we prove that for each
$p\in \mathcal{S}^*\cup\mathcal{S}(\{f_k''\})$, $f_k'$ also converges in
a neighborhood of $p$.

Let $g_{f_k'}=e^{2u_k'}g_{euc}$.
Since $f_k'\in W^{2,2}_{conf}
(D_{4r}(p))$ with $\int_{D_{4r}(p)}|A_{f_k'}|^2d\mu_{f_k'}<8\pi-\tau$ when $r$ is
sufficiently small and $k$ sufficiently large,
by the arguments in \cite{K-L},
we can find a $v_k$ solving the equation
$$-\Delta v_k=K_{f_k'}e^{2u_k'},\s z\in D_r\s and\s \|v_k\|_{L^\infty(D_r(p))}<C.$$
Since $f_k'$ converges to a conformal
immersion in $D_{4r}\setminus D_{\frac{1}{4}r}(p)$,
we may assume that
$$\|u_k'\|_{L^\infty(D_{2r}\setminus
D_r(p))}<C.$$
 Then
$u_k'-v_k$ is a harmonic function with
$\|u_k'-v_k\|_{L^\infty(\partial D_{2r}(p))}<C$,
then we get $\|u_k'(z)-v_k(z)\|_{L^\infty(D_{2r}(p))}<C$
from the Maximum Principle. Thus, $\|u_k'\|_{L^\infty(D_{2r}(p))}<C$,
which implies $\|\nabla f_k'\|_{L^\infty(D_{2r})}<C$.
By the equation $\Delta f_k'=e^{2u_k'}H_{f_k'}$, and
the fact that
$$\|e^{2u_k'}H_{f_k'}\|_{L^2
(D_{2r})}^2<
e^{2\|u_k'\|_{L^\infty(D_{2r}(p))}}\int_{D_{2r}}|H_{f_k'}|^2d\mu_{{f_k'}},$$
we get $\|\nabla{f_k'}\|_{W^{1,2}(D_{r})}<C$.
We complete the proof.

\endproof

{}


\begin{thebibliography}{2}


\bibitem{R-Y} Y. Bernard and T. Rivi\'ere: Energy quantization for Willmore surfaces and applications. {\em Ann. of Math. (2)}
    {\bf 180}  (2014),   87-136.

\bibitem{Br} R. Bryant: A duality theorem for Willmore surfaces,
{\em J. Differential Geom.}, {\bf 20}
(1984), 23-53.

\bibitem{C} B. Y. Chen: Some conformal invariants
of submanifolds and their applications, {\em
Boll. Un.
Mat.  Ital.} {\bf 10} (1974), 380--385.

\bibitem{C-L}J. Chen and T. Lamm:
A Bernstein type theorem for entire Willmore graphs.
{\em J. Geom. Anal.}  {\bf 23}  (2013), 456-469.


\bibitem{D-K}D. DeTurck and J. Kazdan:
Some regularity theorems in Riemannian geometry.
{\em Ann. Sci. cole Norm. Sup. (4)} {\bf 14}
 (1981),  249--260.


\bibitem{H}F. H\'elein: Harmonic maps,
conservation laws and moving frames.
Translated from the 1996 French original.
With a foreword by James Eells. Second edition.
Cambridge Tracts in Mathematics, 150.
Cambridge University Press, Cambridge, 2002.

\bibitem{Hu}A. Huber: On subharmonic functions and
differential geometry in the large,
{\em Comment. Math. Helv.} {\bf 32} (1957), 181-206.



\bibitem{K-L} E. Kuwert and Y. Li: $W^{2,2}$-conformal immersions of a closed Riemann
surface into $\R^n$.  {\em Comm. Anal. Geom.} {\bf 20} (2012),  313-340.

\bibitem{K-S}
E. Kuwert and R. Sch\"atzle:
The Willmore flow with small initial energy.
{\em J. Differential Geom.}, {\bf 57} (2001), 409-441.



\bibitem{K-S2} E. Kuwert and R. Sch\"atzle: Removability
of point singularities of Willmore surfaces, {\em Ann. of Math.}
{\bf 160} (2004), 315-357.

\bibitem{L-Y} P. Li and S.T. Yau: A new conformal
invariant and its applications to the Willmore conjecture
and the first eigenvalue on compact surfaces,
{\em Invent. Math} {\bf 69} (1982), 269-291.


\bibitem{L}Y. Li: Weak limit of an immersed surface sequence with bounded Willmore functional. {\em arXiv:1109.1472}.

\bibitem{L-S} Y. Luo and J. Sun:
Remarks on a Bernstein type theorem for entire
Willmore graphs in $R^3$.
{\em J. Geom. Anal.}  {\bf 24}  (2014),  1613-1618.

\bibitem{M-S} S. M\"uller and V. \v{S}ver\'ak: On surfaces
of finite total curvature, {\em J. Differential Geom.}
{\bf 42} (1995),
229-258.

\bibitem{P-R} J. P\'erez and A. Ros:
Properly embedded minimal surfaces with finite total curvature.
The global theory of minimal surfaces in flat spaces
(Martina Franca, 1999), 15-66,
{\em Lecture Notes in Math.}, {\bf 1775}, Springer, Berlin, 2002.


\bibitem{R} T. Rivi\'ere:  Analysis aspects of Willmore
surfaces. {\em Invent. Math.}  {\bf 174 } (2008),  no. 1, 1-45.

\bibitem{R2} T. Rivi\'ere: Lipschitz conformal immersions from degenerating Riemann surfaces with $L^2$-bounded second fundamental forms. {\em Adv. Calc. Var.}  {\bf 6}  (2013), 1-31.

\bibitem{S} L. Simon: Existence of surfaces minimizing
the Willmore functional, {\em Comm. Anal. Geom.}
{\bf 1} (1993),
281-326.


\bibitem{W} T. J. Willmore: Total Curvature in Riemannian
Geometry, John Wiley \& Sons, New York (1982).



\end{thebibliography}
\end{document}